\documentclass[a4paper,11pt]{article}

\usepackage{amsmath, amssymb, amsfonts, amsthm, diagrams}
\usepackage[latin1]{inputenc}
\usepackage[english]{babel}

\theoremstyle{plain}
\newtheorem{teo}{Theorem}[section]
\newtheorem{lem}[teo]{Lemma}
\newtheorem{prop}[teo]{Proposition}
\newtheorem{cor}[teo]{Corollary}

\theoremstyle{definition}

\newtheorem{rem}[teo]{Remark}

\newcommand{\sln}{{\mathfrak{sl}(n)}}
\newcommand{\gln}{{\mathfrak{gl}(n)}}

\newcommand{\Cbb}{{\mathbb C}}
\newcommand{\Qbb}{{\mathbb Q}}
\newcommand{\Zbb}{{\mathbb Z}}
\newcommand{\Nbb}{{\mathbb N}}
\newcommand{\Ocal}{{\mathcal O}}
\newcommand{\Rcal}{{\mathcal R}}

\newcommand{\psr}{{      \psi_{\rm  reg}        }}
\newcommand{\pr}{{      \phi_{\rm  reg}        }}
\newcommand{\gr}{{      \mathfrak{g}_{\rm reg}        }}
\newcommand{\hr}{{      h_{\rm reg}        }}
\newcommand{\Hr}{{\operatorname{H}_ {  \rm reg      } }}
\newcommand{\Xr}{{\mathfrak{X}_{\rm  reg} }}
\newcommand{\Ns}{{\mathcal S_\sigma}}
\newcommand{\Nss}{{\mathcal S}}
\newcommand{\Hsp}{{\operatorname{H}_\sigma}}
\newcommand{\Hs}{{\operatorname{H}_\sigma^\prime}}
\newcommand{\Zs}{{\operatorname{Z}_\sigma}}

\newcommand{\ab}{{ {\mathbb A} ^{b_1-1}  }}
\newcommand{\ds}{{  /\hspace{-.5ex}/          }}
\newcommand{\Xcal}{{\mathfrak{X}}}

\newcommand{\Pbb}{{\mathbb P}}
\newcommand{\Abb}{{\mathbb A}}

\newcommand{\Gbb}{{\mathbb G _m}}
\newcommand{\Fcal}{{\mathcal F}}

\author{S{\'e}bastien Jansou and Nicolas Ressayre}
\title{INVARIANT DEFORMATIONS\\ OF ORBIT CLOSURES IN $\sln$}

\begin{document}
\date{}
\maketitle

\begin{abstract}

We study deformations of orbit closures for the action of a connected semisimple group $G$ on its Lie algebra $\mathfrak{g}$, especially when $G$ is the special linear group.

The tools we use are on the one hand the invariant Hilbert scheme and on the other hand the sheets of $\mathfrak{g}$.
We show that when $G$ is the special linear group, the connected components of the invariant Hilbert schemes we get are the geometric quotients of the sheets of $\mathfrak{g}$. These quotients were constructed by Katsylo for a general semisimple Lie algebra $\mathfrak{g}$; in our case, they happen to be affine spaces.

\end{abstract}

\section*{Introduction}

Let $G$ be a complex reductive group, and $V$ be a finite dimensional $G$-module. A fundamental problem is to endow some sets of orbits of $G$ in $V$ with a structure of variety. The geometric invariant theory is the classical answer in this context: the set of closed orbits of $G$ in $V$ has a natural structure of affine variety. We denote by $V\ds G$ this variety, equipped with a $G$-invariant quotient map $\pi : V \rightarrow V\ds G$.

Recently, Alexeev and Brion defined in \cite{AB} a structure of quasiprojective scheme on some sets of $G$-stable closed affine subscheme of $V$. A natural question is to wonder what happens when one applies Alexeev-Brion's construction to the orbit closures of $G$ in $V$. Here, we study this construction in the case of a well known $G$-module, namely the adjoint representation of a semisimple group $G$, especially when $G$ is the special linear group $\operatorname{SL}(n)$.

From now on, we assume that $G$ is semisimple, and denote by $\mathfrak{g}$ its Lie algebra endowed with the adjoint action of $G$. 
Let us recall that a sheet of $\mathfrak{g}$ is an irreducible component of the set of points in $\mathfrak{g}$ whose $G$-orbit has a fixed dimension. Let us fix a sheet $\Nss$. We show that the $G$-module structure on the affine algebra $\Cbb[\overline{G \cdot x}]$ of the orbit closure $\overline{G \cdot x}$ of $x$ doesn't depend on $x$ in $\Nss$. This allows us to define a set-theoretical application from $\Nss$ to some Alexeev-Brion's invariant Hilbert scheme of $\mathfrak{g}$: $$\begin{array}{cclc} \pi_\Nss : & \Nss & \longrightarrow&  \operatorname{Hilb}^G_\Nss(\mathfrak{g}) \\ ~ & x &\longmapsto & \overline{G \cdot x}.
\end{array}$$
A unique sheet is open in $\mathfrak{g}$: we call it the \textit{regular} one, and denote it by $\gr$. 

In Section~\ref{sec:greg} we are interested in $\operatorname{Hilb}^G_{\gr}(\mathfrak{g})$.
The graph of the quotient map $\pi\,:\,\mathfrak{g}\rightarrow\mathfrak{g}\ds G$ is a flat family of 
$G$-stable closed subschemes of $\mathfrak{g}$ over $\mathfrak{g}\ds G$.
So, this family is the pullback of the universal one by a morphism. We prove that this morphism 
is an isomorphism by showing that $\operatorname{Hilb}^G_{\gr}(\mathfrak{g})$ is smooth and applying Zariski's main theorem. 
So, we obtain that the application $\pi_{\gr}$ identifies with the restriction of the quotient map 
$\pi:\mathfrak{g}\rightarrow \mathfrak{g}\ds G$; in particular, it is a morphism.

In Section~\ref{sec:sln}, we study any sheet $\Nss$ for $G=\operatorname{SL}(n)$. 
We explicitly construct a flat family over an affine space whose fibers are the closures 
in $\mathfrak{g}$ of the $G$-orbits in $\Nss$.
Then, we show following the same method as in the case of $\gr$ that this family is universal. 
Let us denote by $\pi:\Nss \rightarrow \Nss/\operatorname{SL}(n)$ the geometric quotient of $\Nss$, constructed by Katsylo in \cite{K}. We show that there is a canonical morphism

 $$
\begin{array}{cclc} \theta : & \Nss / \operatorname{SL}(n) & \longrightarrow &  \operatorname{Hilb}^{\operatorname{SL}(n)}_\Nss(\mathfrak{g}) \\ 
~ &  \operatorname{SL}(n) \cdot x &\longmapsto & \overline{ \operatorname{SL}(n) \cdot x} 
\end{array}$$ which is actually an isomorphism onto a connected component of $\operatorname{Hilb}^{\operatorname{SL}(n)}_\Nss(\mathfrak{g})$.

Another motivation for this work is to understand examples of invariant Hilbert schemes. Indeed, the construction of Alexeev and Brion is indirect and only few examples are known (see \cite{J}, \cite{BC}). Here, the connected components of invariant Hilbert schemes we obtain happen to be affine spaces, as in \cite{J} and \cite{BC}. Note that this answers in the case of $\operatorname{SL}(n)$ to a question of Katsylo who asked if the geometric quotient $\Nss/G$ is normal.

\section{Hilbert's sheets}

We consider schemes and affine algebraic groups over $\Cbb$.
Let $G$ be a connected semisimple group. We choose a Borel subgroup $B$, and a maximal torus $T$ contained in $B$. We denote by $U$ the unipotent radical of $B$; we have $B=TU$.

We denote by $\Lambda$ the character group of $T$. 
We denote by $\Lambda^+$ the set of elements of $\Lambda$ that are dominant weights with respect to $B$.
The set $\Lambda ^ +$ is in bijection with the set of isomorphism classes of simple rational $G$-modules. If $\lambda$ is an element of $\Lambda^+$, we denote by $V(\lambda)$ a simple $G$-module associated, that is of highest weight $\lambda$.

If $V$ is a rational $G$-module, we denote by $V_{(\lambda)}$ its isotypical component of type $ \lambda$, that is the sum of its submodules isomorphic to $ V(\lambda ) $. We have the decomposition $V=\bigoplus_{ \lambda \in  \Lambda ^ +} V_{(\lambda)} $.

In any decomposition of $V$ as a direct sum of simple modules, the multiplicity of the simple module  $V(\lambda)$ is the dimension of $V_{(\lambda)}^U$. We say that $V$ \textit{has finite multiplicities} if these multiplicities are finite (for any dominant weight $\lambda$).\\

Let us recall some definitions from \cite[\S1]{AB}. A \textit{family of affine $G$-schemes} over some scheme $S$ is a scheme $ \mathfrak{X}$ equipped with an action of $G$ and with a morphism $\pi:\mathfrak{X} \rightarrow S$ that is affine, of finite type and $G$-invariant. We have a $G$-equivariant morphism of $\Ocal _ S$-modules $$ \pi _ * \Ocal _ { \mathfrak{X} } \simeq \bigoplus _ {\lambda \in \Lambda^+} \Fcal _ \lambda \otimes_\Cbb V ( \lambda),$$  where each $\Fcal _ \lambda :=     (\pi _ * \Ocal _ { \mathfrak{X} } )_{(\lambda)}^U$ is equipped with the trivial action of $G$. Let $h :\Lambda ^ + \rightarrow \Nbb$ be a function. The family $\mathfrak{X}$ is said to be \textit{of Hilbert function $h$} if  each  $\Fcal _ \lambda$       is an $\Ocal _ S$-module locally free of rank $h(\lambda)$. (Then the morphism $\pi$ is flat.)


Let $X$ be an affine $G$-scheme, and $h :\Lambda ^ + \longrightarrow \Nbb$ a function.
A \textit{family of $G$-stable closed subschemes of $X$} over some scheme $S$ is a  $G$-stable closed subscheme $ \mathfrak{X} \subseteq S \times X$. The projection $S \times X \rightarrow S$ induces a family of affine $G$-schemes $\mathfrak{X} \rightarrow S $.
The contravariant functor: $( \mbox{Schemes}) ^ \circ \longrightarrow (\mbox{Sets})$ that associates to every scheme $S$ the set of families $ \mathfrak{X} \subseteq  S \times X$ of Hilbert function $h$ is  represented by a quasiprojective scheme denoted by $Hilb_h^G(X)$ (\cite[\S1.2]{AB}.\\

The dimension of an affine $G$-scheme whose affine algebra has finite multiplicities can be read on its Hilbert function:
\begin{prop} \label{dim} Let $h :\Lambda ^ + \longrightarrow \Nbb$ be a function. Let $Y$ and $Z$ be two affine schemes of Hilbert function $h$. Then $\dim Y = \dim Z$.
\end{prop}

\begin{proof} Let us denote by $A$ the affine ring of $Y$. 

If $Y$ is \textit{horospherical}, that is (\cite[Lemma 2.4]{AB}) if for any dominant weights $\lambda$, $\mu$, we have $A_{(\lambda)} \cdot A_{(\mu)} \subseteq A_{(\lambda+\mu)}$, it is clear that the dimension of $Y$ can be read on its Hilbert function. Indeed, let us denote by $\theta_0$ the linear map from $\Lambda \otimes \Qbb$ to $\Qbb$ which associates to any fundamental weight the value $1$. We denote by $\theta$ the group homomorphism from $\Lambda$ to $\Zbb$ that is the restriction of $\theta_0$. We associate to $\theta$ a graduation of the algebra $A$ by $\Nbb$: its homogeneous component of degree $d$ is $$A_d:=\bigoplus_{\lambda \in \Lambda^+,~\theta(\lambda) =d} A_{(\lambda)}.$$ The dimension of $A_d$ is finite, and depends only on $h$: $$\dim A_d = \sum_{\lambda \in \Lambda^+,~\theta(\lambda) =d}h(\lambda) \dim V(\lambda).$$ So the Hilbert polynomial of the graded algebra $A$ depends only on $h$, and so does the dimension of $Y$.

We can deduce the proposition. Indeed, $Y$ admits a flat degeneration over a connected scheme to a horospherical $G$-scheme $Y^\prime$ that admits the same Hilbert function (by \cite[Theorem 2.7]{AB}). So $\dim Y = \dim Y ^\prime$ depends only on $h$. \end{proof}

We will use the method of ``asymptotic cones'' of Borho and Kraft (\cite[\S5.2]{PV}): let $V$ be a finite dimensional rational $G$-module and $F$ the closure of an orbit in $V$ (or, more generally, any $G$-stable closed subvariety contained in a fiber of the categorical quotient $V \rightarrow V\ds G$). We embed $V$ into the projective space $\Pbb(\Cbb \oplus V)$ of vector lines of $\Cbb \oplus V$ by the inclusion $v \mapsto [1 \oplus v].$ The closure of $F$ in  $\Pbb(\Cbb \oplus V)$ is denoted by $\overline{F}$. The affine cone in  $\Cbb \oplus V$ over $\overline{F}$ is the closed cone $ \mathfrak{X}$ generated by $F$.

The vector space $\Cbb \oplus V$, equipped with its natural scheme structure, is denoted by $\Abb^1 \times V$. The cone $\mathfrak{X} \subseteq \Abb^1 \times V$, viewed as a reduced closed subscheme, is a flat family of affine $G$-schemes. Its fibers over non-zero elements are homothetic to $F$. Its fiber over $0$ is a reduced cone, denoted by $\hat{F}$. It is contained in the null-cone of $V$ (that is the fiber of the categorical quotient $V \rightarrow V\ds G$ containing $0$). Its dimension is the same as $F$.\\

We consider the adjoint action of $G$ on its Lie algebra $\mathfrak{g}$. If $x$ is an element of $\mathfrak{g}$, the affine algebra of the closure of its orbit, viewed as a reduced scheme, has finite multiplicities. Let us denote by $h_x$ its Hilbert function; we call it the Hilbert function associated to $x$. In this paper, we are interested in the connected component denoted $\operatorname{Hilb}^G_x$ of the scheme $\operatorname{Hilb}_{h_x}^G(\mathfrak{g})$ that contains $\overline{G \cdot x}$. It gives the $G$-invariant deformations of $\overline{G \cdot x}$ embedded in $\mathfrak{g}$. We determine it when $x$ is in $\gr$ in \S2, and for any $x$ when $G$ is the special linear group in \S3.

Let us denote by $G_x$ the stabilizer of $x$ in $G$, and $\mathfrak{g}_x$ its Lie algebra. The coadjoint action of $G_x$ is its natural action on the dual vector space $\mathfrak{g}_x^*$.

\begin{prop} \label{et} Let us assume the orbit closure $\overline{G \cdot x}$ to be normal. The tangent space $T_{\overline{G \cdot x}} \operatorname{Hilb}^G_x$ to $\operatorname{Hilb}^G_x$ at the point  $\overline{G \cdot x}$ is canonically isomorphic to the space of invariants of the coadjoint action of $G_x$.
\end{prop} 

\begin{proof}

The tangent space to $\overline{G \cdot x}$ at the point $x$ is $\mathfrak{g}.x$; it is stable under the action of $G_x$. We denote by 
$[\mathfrak{g} / \mathfrak{g}.x]^{G_{x}}$ the space of invariants under the action of $G_x$ on the quotient vector space $\mathfrak{g} / \mathfrak{g}.x$. According to \cite[Proposition 1.15 (iii)]{AB}, we have a canonical isomorphism    \begin{equation} \label{2}   T_{\overline{G \cdot x}} \operatorname{Hilb}^G_x \cong [\mathfrak{g} / \mathfrak{g}.x]^{G_{x}}. \end{equation} Indeed, the orbit closure $\overline{G \cdot x}$ is assumed to be normal. Moreover, every orbit in $\mathfrak{g}$ has even dimension, and has a finite number of orbits in its closure (\cite[Corollary 3 page 198]{PV}), so the codimension of the boundary of $G \cdot x$ in $\overline{G \cdot x}$ is at least $2$, and the proposition of \cite{AB} can be applied.

To transform (\ref{2}) into the isomorphism of the proposition, we will use the Killing form on $\mathfrak{g}$, denoted by $\kappa$. As $\mathfrak{g}$ is semisimple, its Killing form gives an isomorphism  $$\begin{array}{cclc} \phi : & \mathfrak{g} & \longrightarrow&  \mathfrak{g}^* \\ ~ & y &\longmapsto & \kappa(y,\cdot).
\end{array}$$
The isomorphism $\phi$ is $G$-equivariant, thus $G_x$-equivariant. It sends $\mathfrak{g}.x$ onto the space $\mathfrak{g}_x^\perp$ of linear forms on $\mathfrak{g}$ that vanish on $\mathfrak{g}_x$. Indeed, the common zeros of the elements of $\phi(\mathfrak{g}.x)$ are the elements $y$ in $\mathfrak{g}$ such that $$\forall z \in \mathfrak{g},~ \kappa([z,x],y)=0,$$ that is $$\forall z \in \mathfrak{g},~ \kappa(z,[x,y])=0,$$ and this last condition means that $y$ belongs to $\mathfrak{g}_x$ since $\kappa$ is non-degenerate.

Thus the short exact sequence of $G_x$-modules $$0 \longrightarrow \mathfrak{g}.x \longrightarrow \mathfrak{g} \longrightarrow  \mathfrak{g}/ \mathfrak{g}.x \longrightarrow 0$$ identifies (thanks to $\phi$) with  $$0 \longrightarrow \mathfrak{g}_x^\perp \longrightarrow \mathfrak{g}^* \longrightarrow  (\mathfrak{g}_x)^* \longrightarrow 0,$$ and the proposition follows from (\ref{2}).
\end{proof}

A \textit{sheet} of $\mathfrak{g}$ is a maximal irreducible subset of $\mathfrak{g}$ consisting of $G$-orbits of a fixed dimension. Every sheet of $\mathfrak{g}$ contains a unique nilpotent orbit. A \textit{regular element} of $\mathfrak{g}$ is an element of  $\mathfrak{g}$ whose orbit has maximal dimension. The open subset of  $\mathfrak{g}$ whose elements are the regular elements is a sheet denoted by $\gr$.

Let us call \textit{Hilbert's sheet} a maximal irreducible subset of $\mathfrak{g}$ consisting of elements admitting a fixed associated Hilbert function.

\begin{prop} \label{nh} The Hilbert's sheets of $\mathfrak{g}$ coincide with its sheets.
\end{prop}

\begin{proof} According to Proposition \ref{dim}, any Hilbert's sheet is contained in some sheet. It just remains to check that two points of some sheet $\Nss$ have the same associated Hilbert function.

Let $F$ be the closure of an orbit in $\Nss$. We recalled that its asymptotic cone $\hat{F}$ is a degeneration of $F$. In particular, it is contained in the closure of $\Nss$. Moreover, $\hat{F}$ is contained in the null-cone of $\mathfrak{g}$, and its dimension is the same as $F$. So $\hat{F}$ is the closure of the nilpotent orbit of $\Nss$.

The affine algebra of $\mathfrak{g}$ is the symmetric algebra of $\mathfrak{g}^*$. Its graduation induces a $G$-invariant filtration on the affine algebra $A$ of $F$. The affine algebra of the asymptotic cone $\hat{F}$ is isomorphic, as an algebra equipped with an action of $G$, to the graded algebra $\hat{A}$ associated to the filtered algebra $A$. In particular, $A$ and $\hat{A}$ are isomorphic as $G$-modules, and their multiplicities are equal: the Hilbert function of $F$ is equal to that of $\hat{F}$, and the proposition is proved. \end{proof}

Notice that in the case of the regular sheet, Proposition~\ref{nh} is a direct consequence of \cite[Theorem 0.9]{Ko}.

\section{Regular case}
\label{sec:greg}

Let us denote by $\hr$ the Hilbert function associated to the regular elements of $\mathfrak{g}$ (Proposition \ref{nh}). In this section, we prove that the invariant Hilbert scheme $\Hr:=\operatorname{Hilb}_\hr^G(\mathfrak{g})$ is the categorical quotient $\mathfrak{g}\ds G$, that is an affine space whose dimension is the rank of $G$.

By \cite[Theorem 0.1]{Ko}, all schematic fibers of the quotient morphism $\mathfrak{g}\rightarrow\mathfrak{g}\ds G$ are reduced.
This allows us to  identify in the following the schematic fibers with the set-theoretical fibers.

\subsection{A morphism from $\mathfrak{g}\ds G$ to $\Hr$}

Let $\Xr$ be the graph of the canonical projection $\mathfrak{g} \rightarrow \mathfrak{g}\ds G$. It is a family of $G$-stable closed subschemes of $\mathfrak{g}$ over $\mathfrak{g}\ds G$.

\begin{prop} \label{morph0} The closed subscheme $\Xr$ is a family of $G$-stable closed subschemes of $\mathfrak{g}$ with Hilbert function $\hr$.
\end{prop}

\begin{proof} Let us denote by $\pi : \Xr \rightarrow \mathfrak{g} \ds G$ the canonical projection, and by $\Rcal:=\pi_* \Ocal _ \Xr$ the direct image by $\pi$ of the structural sheaf of $\Xr$. We have to prove that for any dominant weight $\lambda$, we have that $\Rcal_{(\lambda)}^U$ is a locally free sheaf on $\mathfrak{g} \ds G$ of rank $h(\lambda)$.

Let us first study the case where $\lambda =0$. The morphism $\pi \ds G : \Xr \ds G \rightarrow \mathfrak{g} \ds G$ induced by $\pi$ is clearly an isomorphism. So $\Rcal^G=\Rcal_{(0)}^U$ is a free module on $\mathfrak{g} \ds G$ of rank $1=\hr(0)$.

Let $\lambda$ be a dominant weight. It is known (see \cite[Lemma 1.2]{AB}) that $\Rcal_{(\lambda)}^U$ is a coherent $\Rcal^G$-module. Thus it is a coherent module on $\mathfrak{g} \ds G$. To see that it is locally free, we just have to check that its rank is constant. The fibers of $\pi$ are those of the canonical projection $\mathfrak{g} \rightarrow \mathfrak{g}\ds G$, so they are the orbit closures of the regular elements, and all of them admit $\hr$ as Hilbert function. So the rank of $\Rcal_{(\lambda)}^U$ at any closed point of $\mathfrak{g}\ds G$ is $h(\lambda)$, and the proposition is proved.    \end{proof}

This gives us a canonical morphism $$\pr : \mathfrak{g}\ds G \longrightarrow \Hr.$$
We will prove in the following of \S2 that $\pr$ is an isomorphism.

\begin{lem} \label{bij0} The morphism $\pr$ realizes a bijection from the set of closed points of $\mathfrak{g}\ds G$ to the set of closed points of $\Hr$.
\end{lem}

\begin{proof} We remark that $\pr$ is injective. Let us check it is surjective: in other words, that any $G$-invariant closed subscheme of $\mathfrak{g}$ of Hilbert function $\hr$ is a fiber of $\mathfrak{g} \rightarrow \mathfrak{g}\ds G$.

Let $Y$ be such a subscheme. As $\hr(0)=1$, it has to be contained in some fiber $F$ of $\mathfrak{g} \rightarrow \mathfrak{g}\ds G$ over a reduced closed point. But $F$ already corresponds to a closed point of $\Hr$ in the image of $\pr$. Moreover, $F$ admits no proper closed subscheme admitting the same Hilbert function, so $F=Y$, and the lemma is proved.
\end{proof}

Let us denote by $r$ the rank of $G$. The quotient $\mathfrak{g}\ds G$ is an affine space of dimension $r$. A consequence of Lemma \ref{bij0} is:

\begin{cor} \label{coro7} The dimension of $\Hr$ is $r$.
\end{cor}

\subsection{Tangent space}

In this section, we prove:

\begin{prop} The scheme $\Hr$ is smooth.
\end{prop}

\begin{proof} Let $Z$ be a closed point of $\Hr$. We have to prove that the dimension of the tangent space $T_{Z}\Hr$ is $r$. We still denote by $Z$ the closed subscheme of $\mathfrak{g}$ corresponding to $Z$. By Lemma \ref{bij0}, we know that $Z$ is a fiber of the morphism $\mathfrak{g} \rightarrow \mathfrak{g}\ds G$, thus the closure of some regular element $x$. It is a normal variety. By Proposition \ref{et}, we have to prove that the dimension of $$(\mathfrak{g}_x^*)^{G_x}$$ is $r$, or simply that it is lower or equal to $r$ (by Corollary \ref{coro7}). 

Let us prove that the dimension of the bigger space $$(\mathfrak{g}_x^*)^{\mathfrak{g}_x}$$ is $r$, and the proposition will be proved.

A linear form on $\mathfrak{g}_x$ is $\mathfrak{g}_x$-invariant iff it vanishes on the derived algebra $[ \mathfrak{g}_x , \mathfrak{g}_x ]$, so we have to prove that $$(\mathfrak{g}_x / [ \mathfrak{g}_x , \mathfrak{g}_x ])^*$$ is $r$-dimensional. We will prove that $\mathfrak{g}_x$ is an $r$-dimensional abelian algebra, and the proposition will be proved. This is true if $x$ is semisimple, because then $\mathfrak{g}_x$ is a Cartan subalgebra of $\mathfrak{g}$. If the regular element $x$ is not assumed to be semisimple, the dimension of $\mathfrak{g}_x$ is still $r$, because this doesn't depend on the regular element $x$, by definition. Let us check that $\mathfrak{g}_x$ is abelian.

Let us denote by $\operatorname{Grass}_r(\mathfrak{g})$ the grassmannian of $r$-dimensional subspaces of $\mathfrak{g}$, endowed with its projective variety structure. The subset of $\gr \times \operatorname{Grass}_r(\mathfrak{g})$: $$\{(z,\mathfrak{h}) \in \gr \times \operatorname{Grass}_r(\mathfrak{g}) ~|~ \mathfrak{h} \cdot z=0 \mbox{ and } [\mathfrak{h},\mathfrak{h}]=0 \}$$ is closed, so its image by the natural projection into $\gr$ is closed too. As its image contains the semisimple elements of $\gr$, it is equal to $\gr$. Thus $\mathfrak{g}_x$ is abelian for any regular $x$, and the proposition is proved.    \end{proof}

\subsection{Conclusion}

We can now conclude that the family $\Xr$ of Proposition \ref{morph0} is the universal family:

\begin{teo} \label{theor} The morphism $\pr$ from  $\mathfrak{g}\ds G$ to $\Hr$ is an isomorphism.
\end{teo}

\begin{proof} The morphism $\pr$ is bijective (Lemma \ref{bij0}) and $\Hr$ is normal. According to Zariski's main theorem, $\pr$ is an isomorphism.
\end{proof}

\begin{rem} One knows there is a canonical morphism $$\psr : \Hr \longrightarrow \mathfrak{g} \ds G$$ that associates to any closed point $F$ of $\Hr$ (viewed as a closed subscheme of $\mathfrak{g}$) its categorical quotient $F \ds G$ (viewed as a closed point of $\mathfrak{g} \ds G$). This morphism is a particular case of morphism $$\eta : \operatorname{Hilb}^G_{h}(V) \longrightarrow \operatorname{Hilb}_{h(0)}(V\ds G)$$ defined in \cite[\S1.2]{AB}, because $\hr(0)=1$ and thus the punctual Hilbert scheme that parametrizes closed subschemes of length $1$ in $\mathfrak{g} \ds G$ identifies with $\mathfrak{g} \ds G$ itself. The morphism $\psr$ is clearly the inverse morphism of $\pr$. \end{rem}

\begin{rem} 
As pointed to us by M. Brion, Theorem \ref{theor} admits the following generalization:\\

\textit{Let $X$ be an irreducible affine $G$-variety such that $\pi: X \rightarrow X \ds G$ is flat. Let $h$ be the Hilbert function of its fibers. Then the graph $\Gamma$ of $\pi$ is the universal family; in particular, $\operatorname{Hilb}^G_h(X)$ identifies with $X//G$.}\\

The idea of his proof is to check that $\Gamma$ represents the functor. 
Let $\Xcal \subseteq X \times S$ be a flat family of Hilbert function $h$, over some affine scheme $S$. 
Since $h(0)=1$, the scheme $S$ identifies with $\Xcal \ds G$ and maps on $X \ds G$ 
(by the morphism induced by the first projection $X \times S \rightarrow X$). 
We obtain the following commutative diagram:

\begin{diagram}
 \Xcal&\rTo^{p_2}&\Gamma\\
 \dTo&~&\dTo\\
\Xcal\ds G\simeq S&\rTo&X\ds G.\\
\end{diagram}

It remains to prove that $\Xcal$ is isomorphic (canonically) to the fiber product $\Gamma \times_{X \ds G} S$. 
This has only to be verified over the closed points of $S$.
The assertion follows.\\

The Hilbert schemes we obtain applying the above Brion's result to  $G$-modules are 
always affine spaces.
The representations $V$ of a simple group $G$ such that $V\rightarrow V\ds G$ is flat 
have been classified by G.~Schwarz in \cite{Sch}.\\

Unfortunately, the sheets of $\sln$  are not affine in general and Katsylo's quotient
cannot be extended to their closure.
So, Brion's theorem cannot be applied, 
whereas the method we used to prove Theorem~\ref{theor} can be used.
\end{rem}

\section{Case of $\sln$}
\label{sec:sln}

We denote by $t$ an indeterminate over $\Cbb$, and $I_n$ the identity matrix of size $n \times n$.
If $x$ is an element of $\sln$ and $i=1 \cdots n$, we denote by $Q_i^x(t)$ the monic greatest common divisor (in the ring $\Cbb[t]$) of the $(n+1-i) \times (n+1-i)$-sized  minors of $x-t I_n$, and $Q_{n+1}^x(t):=1$.

Then we put $$q_i^x(t) := Q_i^x(t)/Q_{i+1}^x(t).$$ 
The polynomials $q_1^x(t), \cdots, q_n^x(t)$ are the invariant factors of the matrix $x-tI_n$ with coefficients in the euclidean ring $\Cbb[t]$, ordered in such a way that $q_{i+1}^x(t)$ divides $q_{i}^x(t)$. \\

If $x$, $y$ are elements of $\sln$, then $y$ is in the closure of the orbit $\operatorname{SL}(n) \cdot x $ of $x$ if and only if for any $i=1 \dots n$, the polynomial $Q_i^x(t)$ divides $Q_i^y(t)$. 
In other words, iff for any $i$, the polynomial $Q_i^x(t)$ divides the $(n+1-i) \times (n+1-i)$-sized  minors of $y-t I_n$.

According to \cite{W}, when $x$ is nilpotent, these conditions defines the closure of $\operatorname{SL}(n) \cdot x$ as a reduced scheme: to be more precise, when one divides a  $(n+1-i) \times (n+1-i)$-sized  minor of $y-t I_n$ by $Q_i^x(t)$ using Euclid algorithm, the remainder he gets is a regular function of $y$. All such functions generate the ideal of the closure of $\operatorname{SL}(n) \cdot x$. We will deduce easily from this difficult result that the same remains true if $x$ is no longer assumed to be nilpotent.\\

The set of sheets of $\sln$ is in bijection with the set of partitions $n$, that is   of sequences 
$\sigma=(b_1 \geq b_2 \geq b_3 \geq \dots)$ of nonnegative integers such that $b_1+b_2+b_3+ \dots = n$ (see \cite[\S 2.3]{Bo}).
Namely, if $\sigma$ is a partition of $n$, the elements of the correspondent sheet $\Ns$ are those elements $x$ such that for any $i$, the polynomial $q_i^x(t)$ is of degree $b_i$.
We denote by $\widehat{\sigma}=(c_1 \geq c_2 \geq c_3 \geq \dots)$ the conjugate partition, where $c_j$ is the number of $i$ such that $b_i \geq j$.
We denote by $h_\sigma$ the Hilbert function associated to the points of $\Ns$ (Proposition \ref{nh}).
We denote by $\Zs$ the closure of the nilpotent orbit of $\Ns$. The connected component of $\operatorname{Hilb}^{\operatorname{SL}(n)}_{h_\sigma}(\sln)$ that contains $\Zs$ as a closed point is denoted $\Hsp$. We will prove in this section that $\Hsp$ is an affine space of dimension $b_1-1$. The proof is similar to \S2.

We recall that the sheets of $\sln$ are smooth (\cite{Kr}).

\subsection{A construction of the geometric quotient of $\Ns$}

Katsylo showed in \cite{K} that any sheet of a semisimple Lie algebra admits a geometric quotient. Although his proof contains an explicit construction, it doesn't make clear the geometric properties of the quotient. Here we present a simple description of the quotient in the case of the Lie algebra $\sln$. It takes on the invariant factors theory. We get that the quotient is an affine space.

\begin{lem}
Given some $i$, the application $\Ns \longrightarrow \Abb^{b_i}$ that associates to any $x$ the coefficients of $q_i^x(t)=t^{b_i}+\lambda_{b_i-1}^xt^{b_i-1}+\dots+\lambda_0^x t^0$ is regular.
\end{lem}

\begin{proof}
Up to scalar multiplication, the polynomial $q_i^x(t)$ is the unique nonzero polynomial of degree less or equal to $b_i$ such that \begin{equation} \label{1} \dim \ker q_i^x(x) \geq N:= \sum_{j=1}^{b_i}c_j. \end{equation} Thus the closed subset of $\Ns \times \Pbb^{b_i}$ consisting of elements $(x,[\mu_0:\dots:\mu_{b_i}])$ such that $$\dim \ker (\sum_{j=0}^{b_i}\mu_jx^j) \geq N$$ is the graph of the application  $$\begin{array}{cclc} \psi : & \Ns & \longrightarrow&  \Pbb^{b_i}\\ ~ & x &\longmapsto & [\lambda_0^x:\dots:\lambda_{b_i-1}^x:1]
\end{array}$$

According to \cite[Exercise 7.8 p 76]{Ha}, this graph is also the graph of a rational map $\phi$ from $\Ns$ to $\Pbb^{b_i}$. On the open subset $\Omega$ of $\Ns$ where $\phi$ is regular, $\phi$ coincides with $\psi$, so the functions $x \mapsto \lambda_j^x$ are regular functions from $\Omega$ to $\Abb^1$. As $\Ns$ is smooth, the complementary of $\Omega$ in $\Ns$ has codimension at least 2 (\cite[Thm 3 chap II.3.1]{Sha}). We conclude that the functions extend to regular functions from $\Ns$ to $\Abb^1$. By continuity, these extensions satisfy (\ref{1}), so they coincide with the functions $x \mapsto \lambda_j^x$ on $\Ns$.  \end{proof}

Let us define, for any $x$ in $\Ns$, the monic polynomial of degree $b_i-b_{i+1}$: $$p_i^x(t): = q_i^x(t)/q_{i+1}^x(t)$$ (where $q_{n+1}^x:=1$). It follows from the previous lemma that its coefficients, viewed as functions of $x$, are regular functions from $\Ns$ to $\Abb^1$.

Given an $x$, the family $(p_1^x(t), \dots, p_n^x(t))$ can be any family of monic polynomials of degrees $b_i-b_{i+1}$, provided the following relation is satisfied, where $S(p_i^x)$ denotes the sum of the roots of $p_i^x$, counted with multiplicities (given by its first nondominant coefficient): $$\sum_{i=1}^niS(p_i^x)=0$$ (this relation simply means that the trace of $x$ is zero).

Thus, associating to any $x$ the coefficients of the family $(p_1^x(t), \dots, p_n^x(t))$, we get a regular map $\pi$ from $\Ns$ to a linear hyperplane of $\Cbb^{b_1}$, which we will denote by $\Abb^{b_1-1}$.

\begin{prop} The map $\pi : \Ns \longrightarrow \Abb^{b_1-1}$ is the geometric quotient of $\Ns$.
\end{prop}

\begin{proof} This map is surjective, and its fibers are exactly the orbits of $\Ns$ under the action of $\operatorname{SL}(n)$. Let us denote by $\Ns/\operatorname{SL}(n)$ the geometric quotient of $\Ns$ (whose existence is proved in \cite{K}). The map $\pi$ is the composite of the canonical projection from $\Ns$ to  $\Ns/\operatorname{SL}(n)$ with a regular bijection $$ \Ns/\operatorname{SL}(n) \longrightarrow \Abb^{b_1-1}.$$ This last map is bijective (thus birational), and the space $\Abb^{b_1-1}$ is normal. According to Zariski's main theorem, it is an isomorphism. \end{proof}

\subsection{A morphism from $\Ns/\operatorname{SL}(n)$ to $\Hsp$} \label{paramorph}

If $z=(p_1(t),\dots,p_n(t))$ is a closed point of $\ab$ corresponding to the orbit $\operatorname{SL}(n) \cdot x$ in $\Ns$, the polynomial $$Q_i^x(t)=p_i(t) \cdot (p_{i+1}(t))^2  \cdot ... \cdot (p_{n}(t))^{n-i+1}$$ only depends on $z$. Let us denote it by $Q_i^z(t)$. Its coefficients are regular functions from $\ab$ to $\Abb^1$.

Let us consider the closed subscheme $\Xcal_\sigma$ of $\{(z,y) \in \ab \times \sln \}$ defined by the vanishing, for $i=1\dots n$, of the remainders we get when we divide the $(n+1-i) \times (n+1-i)$-minors of $y-tI_n$ by $Q_i^z(t)$. We denote by $I_\sigma$ the ideal generated by these remainders. The underlying set of $\Xcal_\sigma$ consists of all the couples $(z,y)$ such that $y$ is in the closure of the orbit corresponding to $z$.

\begin{prop} \label{morph} The closed subscheme $\Xcal_\sigma$ is a family of $\operatorname{SL}(n)$-stable closed subschemes of $\sln$ with Hilbert function $h_\sigma$.
\end{prop}

\begin{proof} The proof is similar to that of Proposition \ref{morph0}. The subscheme $\Xcal_\sigma$ is a family of $\operatorname{SL}(n)$-stable closed subschemes of $\sln$ over $\ab$. Let us denote by $\pi$ the morphism $\Xcal_\sigma \longrightarrow \ab$.

As previously, let us first remark that the morphism $$\pi \ds \operatorname{SL}(n) : \Xcal_\sigma \ds \operatorname{SL}(n) \longrightarrow \ab$$ induced by $\pi$ is an isomorphism. To do this, let us verify that the comorphism $$(\pi \ds \operatorname{SL}(n))^* : \Cbb[\ab] \longrightarrow \Cbb[\ab] \otimes \Cbb[\sln]^{\operatorname{SL}(n)} / I_\sigma ^{\operatorname{SL}(n)}$$ is an isomorphism. It is injective, as $\pi$ is surjective. Its surjectivity comes from the relations that define  $\Xcal_\sigma$: they give, for $i=1$, that $Q_1^z(t)$ divides the determinant of $tI_n-y$, that is   the characteristic polynomial of $y$. As their degrees are equal, $Q_1^z(t)$ and the characteristic polynomial of $y$ are equal. This gives the surjectivity.

We go on as previously: let $\lambda$ be a dominant weight. The $R^{\operatorname{SL}(n)}$-module $R_{(\lambda)}^U$ is of finite type (\cite[Lemma 1.2]{AB}). Thus $(\pi_*\Ocal_{\Xcal_\sigma})^U_{(\lambda)}$ is a coherent $\Ocal_{\ab}$-module. To see that it is locally free, we just have to check that its rank is constant. Let us assume that the origin $0 \in \ab$ corresponds to the nilpotent orbit in $\Ns$. The fiber of $\pi$ over $0$ is the closure of this orbit, fitted with its structure of reduced scheme. Thus, the rank of  $(\pi_*\Ocal_{\Xcal_\sigma})^U_{(\lambda)}$ at $0$ is $h_\sigma(\lambda)$. If $z$  is any  point of $\ab$, the fiber of $\pi$ over $z$ is as a set the closure in $\sln$ of the corresponding orbit. So, by Proposition \ref{nh} the rank of $(\pi_*\Ocal_{\Xcal_\sigma})^U_{(\lambda)}$ at $z$ is at least $h_\sigma(\lambda)$. To conclude, we use the action of the multiplicative group on $\sln$ (by homotheties) and the induced action on $\ab$, that makes $\pi$ equivariant. The orbit of $z$ goes arbitrary close to $0$, and the rank of a coherent sheaf is upper semicontinuous, so the rank of $(\pi_*\Ocal_{\Xcal_\sigma})^U_{(\lambda)}$ is $h_\sigma(\lambda)$ at $z$.
\end{proof}


\subsection{Tangent space}

In this section, we compute the dimension of the tangent space to $\Hsp$ at the point $\Zs$:

\begin{prop} The dimension of $T_{\Zs}\Hsp$ is $b_1-1$.
\end{prop}

\begin{proof} Let $x$ be an element in the open orbit in $\Zs$. It is known that $\Zs$ is normal (\cite{KP}). So by Proposition \ref{et}, we just have to prove that the dimension of $$(\sln_x^*)^{\operatorname{SL}(n)_x}$$ is $b_1-1$. Let us consider $\operatorname{SL}(n)$ as a closed subgroup of the general linear group $\operatorname{GL}(n)$, and $\sln$ as a subalgebra of $\gln$. The stabilizer $\operatorname{GL}(n)_x$ of $x$ in $\operatorname{GL}(n)$ is generated by $\operatorname{SL}(n)_x$ and the center of $\operatorname{GL}(n)$. It is clearly equivalent to prove that the dimension of  $$(\gln_x^*)^{\operatorname{GL}(n)_x}$$ is $b_1$. The group $\operatorname{GL}(n)_x$ is connected, so the last space is isomorphic to $$(\gln_x^*)^{\gln_x}.$$ A linear form on $\gln_x$ is $\gln_x$-invariant iff it vanishes on the derived algebra $[ \gln_x , \gln_x ]$, so we have to prove that $$(\gln_x / [ \gln_x , \gln_x ])^*$$ is $b_1$-dimensional. This fact is the following elementary lemma. \end{proof}

\begin{lem} Let $E=\bigoplus_{i=1}^{c_1}E_i$ be a graded vector space over $\Cbb$, where each $E_i$ is $b_i$-dimensional. We denote by $\mathfrak{h}:=\mathfrak{gl}(E)$ the Lie algebra of endomorphisms of $E$. Let $x$ be a nilpotent element of $\mathfrak{h}$ such that each subspace $E_i$ is stabilized by $x$, and the restriction of $x$ to each $E_i$ is cyclic.

Let us denote by $\mathfrak{h}_x$ the stabilizer of $x$ in $\mathfrak{h}$.
Then the codimension of the derived algebra $[\mathfrak{h}_x ,\mathfrak{h}_x ]$ in $\mathfrak{h}_x $ is $b_1$.
\end{lem}

\begin{proof} The graduation of $E$ induces a graduation on the vector space $\mathfrak{h}$: $$\mathfrak{h}=\bigoplus_{i,j} \operatorname{Hom}(E_i,E_j).$$ Let us denote by $p_i:E\longrightarrow E_i$ the natural projections. As they commute with $x$, the subspace $\mathfrak{h}_x$ of $\mathfrak{h}$ is homogeneous: $$\mathfrak{h}_x=\bigoplus_{i,j} \operatorname{Hom}_x(E_i,E_j),$$ where $\operatorname{Hom}_x(E_i,E_j)$ denotes the space of homomorphisms that commute with $x$. Let us choose, for any $i$, an element $e_i$ of $E_i$ such that $x^{b_i-1}e_i \not = 0$. We put $n_{ij}:= b_j-b_i$ if $j<i$ and $0$ otherwise. We denote by $f_{ij}:E_i \rightarrow E_j$ the unique homomorphism that commutes with $x$ and that sends $e_i$ to $x^{n_{ij}}e_j$. Then any homomorphism from $E_i$ to $E_j$ that commutes with $x$ is the composite of $f_{ij}$ with a polynomial in $x$: $$\operatorname{Hom}_x(E_i,E_j)=\Cbb[x] \cdot f_{ij}.$$

We notice that if $i \not = j$, then $\operatorname{Hom}_x(E_i,E_j)$ is contained in $[\mathfrak{h}_x ,\mathfrak{h}_x ]$.
Indeed, for any $u:E_i\rightarrow E_j$, we have $[u,p_i]=u.$

So we have to prove that the codimension in $\bigoplus_{i} \operatorname{Hom}_x(E_i,E_i)$ of $$[\mathfrak{h}_x ,\mathfrak{h}_x ]  \cap \bigoplus_{i} \operatorname{Hom}(E_i,E_i)$$ is $b_1$. The last vector space is generated by its elements of the form $$P(x)[f_{ji},f_{ij}]=P(x)x^{|b_i-b_j|}(\operatorname{id}_{E_i}-\operatorname{id}_{E_j}),$$ where $P(x)$ is a polynomial in $x$.

One checks easily that a basis of a supplementary in $\bigoplus_{i} \operatorname{Hom}_x(E_i,E_i)$ of this space is given by the family of elements $$x^k \operatorname{id}_{E_i}$$ where $0\leq k < b_i-b_{i+1}$, and the lemma is proved.
\end{proof}

\subsection{Conclusion}

In this section, we prove that the family $\Xcal_\sigma$ of Proposition \ref{morph} is the universal family:

\begin{teo} \label{theo7} The morphism $\phi_\sigma$ from  $\Ns/\operatorname{SL}(n)$ to $\Hsp$ obtained in \S\ref{paramorph} is an isomorphism.
\end{teo}

We denote by $\overline{\Ns}$ the closure of $\Ns$ in $\sln$, equipped with its reduced scheme structure. The invariant Hilbert scheme $\Hs := \operatorname{Hilb}^{\operatorname{SL}(n)}_{h_\sigma}(\overline{\Ns})$ which parametrizes the closed subschemes of $\overline{\Ns}$ of Hilbert function $h_\sigma$ is canonically identified with a closed subscheme of $\operatorname{Hilb}^{\operatorname{SL}(n)}_{h_\sigma}(\sln)$. The morphism $\phi_\sigma$ factorizes by a morphism $\psi_\sigma:\Ns/\operatorname{SL}(n) \rightarrow \Hs$. 

To prove the theorem, we will get that the morphism $\psi_\sigma$ is an isomorphism from $\Ns/\operatorname{SL}(n)$ to $\Hs$ and that $\Hs$ is a connected component of $\Hsp$ (Corollary \ref{cordim}).

\begin{lem} \label{bij} The morphism $\psi_\sigma$ induces a bijection from the set of closed points of $\Ns/\operatorname{SL}(n)$ to the set of closed points of $\Hs$.
\end{lem}

\begin{proof} We know that $\psi_\sigma$ is injective. Let us check it is surjective: in other words, that any $\operatorname{SL}(n)$-invariant closed subscheme of $\overline{\Ns}$ with Hilbert function $h_\sigma$ is the closure of some orbit in $\Ns$.

Let $X$ be such a subscheme. As $h_\sigma(0)=1$, it has to be contained in some fiber $F$ of the categorical quotient $\overline{\Ns} \rightarrow \overline{\Ns}\ds \operatorname{SL}(n)$ over a reduced closed point. But $F$ already corresponds to a closed point of $\Hs$ in the image of $\psi_\sigma$. Moreover, $F$ admits no proper closed subscheme admitting the same Hilbert function, so $F=X$, and the lemma is proved.
\end{proof}

\begin{cor} The dimension of $\Hs$ is $b_1-1$.
\end{cor}

The action of the multiplicative group $\Gbb$ on $\sln$ by homotheties induces canonically an action of $\Gbb$ on $\Hsp$, and on $\Hs$ (because it stabilizes $\overline{\Ns}$). The cone $\Zs$ is a $\Gbb$-fixed point of $\Hs$. In fact, it is in the closure of the $\Gbb$-orbit of any point of $\Hs$:

\begin{prop} Let $F$ be a closed point of $\Hs$.
    The morphism $\eta : \Gbb \longrightarrow \Hs$, $t \longmapsto t.X$ extends to a morphism $\Abb ^1 \longrightarrow \Hs$, $0 \longmapsto \Zs$.\end{prop}

\begin{proof} The point $F$ corresponds to a $\operatorname{SL}(n)$-invariant closed subscheme of $\overline{\Ns}$ admitting Hilbert function $h_\sigma$. We still denote it by $F$. As $h_\sigma(0)=1$, it is contained in the fiber of the categorical quotient $\sln \rightarrow \sln\ds {\operatorname{SL}(n)}$ over some closed point. Thus we can apply to it the method of asymptotic cones: we obtain a flat family over $\Abb^1$ whose fiber over $0$ must be $\Zs$ (as in the proof of Proposition \ref{nh}). It gives a morphism from $\Abb^1$ to $\Hs$ whose restriction outside $0$ is $\eta$.
\end{proof}

From the proposition, we deduce that the dimension of the tangent space to $\Hsp$ at any point of $\Hs$ is lower or equal to that at $Z_\sigma$, that is $b_1-1$. As the dimension of $\Hs$ is $b_1-1$, we get:

\begin{cor} \label{cordim}~

 \begin{itemize} \item The scheme $\Hs$ is reduced and smooth.
\item It is a connected component of $\Hsp$.
\end{itemize}
\end{cor} 

The morphism $\psi_\sigma$ is bijective (Lemma \ref{bij}) and $\Hs$ is normal. According to Zariski's main theorem, $\psi_\sigma$ is an isomorphism.
So Theorem \ref{theo7} is proved, thanks to the second point of Corollary \ref{cordim}.


\begin{thebibliography}{999}

\bibitem[AB]{AB} V. Alexeev and M. Brion \\
{\sl Moduli of affine schemes with reductive group action,} J. Algebraic Geom. \textbf{14}, p 83-117, 2005. 

\bibitem[Bo]{Bo} K. Bongartz\\ 
{\sl Schichten von Matrizen sind rationale Varietäten,}  Math. Ann.  \textbf{283},  no. 1,p 53-64, 1989. 

\bibitem[BC]{BC} P. Bravi and S. Cupit-Foutou \\
{\sl Equivariant deformations of the affine multicone over a flag variety,} preprint available on math.AG/0603690. 

\bibitem[Hr]{Ha} J. Harris\\
{\sl Algebraic geometry: a first course}, GTM 133, Springer Verlag, 1992.
 
\bibitem[J]{J} S. Jansou \\
{\sl D{\'e}formations des c{\^o}nes de vecteurs primitifs,}  Math. Ann. \textbf{338}, no 3, p 627-667, 2007. 


\bibitem[Ka]{K} P.I. Katsylo\\
{\sl Sections of sheets in a reductive algebraic Lie algebra,}  Izv. Akad. Nauk SSSR Ser. Mat.   \textbf{46} no. 3, p 477-486, 1982.

\bibitem[Ko]{Ko} B. Kostant\\ 
{\sl Lie group representations on polynomial rings,}  Amer. J. Math.  \textbf{85}, p 327-404, 1963. 

\bibitem[Kr]{Kr}  H. Kraft\\
{\sl Parametrisierung von Konjugationsklassen in $\sln$,}   Math. Ann.   \textbf{234}  no. 3, p 209-220, 1978.

\bibitem[KP]{KP}  H. Kraft and C. Procesi\\
{\sl Closures of conjugacy classes of matrices are normal,} Invent. Math.    \textbf{53}  no. 3, p 227-247, 1979.


\bibitem[PV]{PV} V. Popov and E. Vinberg\\
{\sl Invariant Theory.}\
Encyclopaedia of Mathematical Sciences, vol 55, p 123-278, Springer Verlag 1994.
  
\bibitem[Sch]{Sch} G. W. Schwarz\\
{\sl Representations of simple Lie groups with a free module of covariants,}
Invent. Math.  \textbf{50}  no. 1, p 1-12,  1978/79.

\bibitem[Sha]{Sha} I.R. Shafarevich\\
{\sl Basic algebraic geometry. 1. Varieties in projective space}, Springer Verlag, 1994.
 

\bibitem[W]{W} J. Weyman\\
{\sl The equations of conjugacy classes of nilpotent matrices}, Invent. Math. \textbf{98} no. 2, p 229-245, 1989. 




\end{thebibliography}
\end{document}